\documentclass[12pt]{elsarticle}

\usepackage{commath}
\usepackage{amssymb}
\usepackage{tikz}
\usetikzlibrary{arrows}
\usepackage{hyperref}

\hypersetup{linkcolor=blue,citecolor=red,filecolor=green,urlcolor=magenta} 

\newtheorem{definition}{Definition}[section]
\newtheorem{theorem}{Theorem}[section]
\newtheorem{remark}{Remark}[section]
\newtheorem{lemma}{Lemma}[section]

\journal{Journal of Number Theory}

\begin{document}

\begin{frontmatter}

\title{The Selberg--Delange method and mean value of arithmetic functions over short intervals}

\author{Amrinder Kaur}
\ead{amrinder1kaur@gmail.com}

\author{Ayyadurai Sankaranarayanan}
\ead{sank@uohyd.ac.in}
\affiliation{organization={School of Mathematics and Statistics},
            addressline={University of Hyderabad}, 
            city={Hyderabad},
            postcode={500046}, 
            state={Telangana},
            country={India}}

\date{}

\begin{abstract}
In this paper, we establish a mean value result of arithmetic functions over shorter intervals with the Selberg--Delange method using the Hooley--Huxley contour.
\end{abstract}


\end{frontmatter}

\begin{section}{Introduction} \hfill\\

\noindent
A classical problem in analytic number theory is to study the behaviour of the sum $\sum \limits_{n \leq x} \mu (n)$. It is well known that a bound of the kind 
$$\sum \limits_{n \leq x} \mu (n) \ll x^{\frac{1}{2}+\epsilon} \qquad \text{for every} \; \; \epsilon >0 $$ 
is equivalent to the unproven Riemann hypothesis. \\
More generally, given an arithmetical function $f(n)$, studying the behaviour of the sum $\sum \limits_{n \leq x} f(n)$ is a classical problem. If one knows the analytic properties of the $L$--function attached to $f(n)$, namely 
$$\sum \limits_{n=1}^{\infty} \frac{f(n)}{n^s}$$
(particularly certain growth conditions) and if one knows the nature of the singularity (particularly having only real poles), then Perron's formula \cite{Op} is an appropriate tool to obtain the asymptotic nature of the required sum with a possible good error term. However, if the $L$--function has some singularities whose nature is unknown and having some natural product representation, then Selberg \cite{Ase} and later Delange \cite{Hd1,Hd2} developed a method that enables us to study the sum in question in detail. \\

\noindent
Throughout the paper, the constants $a$ with suffixes are positive constants that need not be the same at each occurrence.

\noindent
In this paper, we consider $\mathcal{P}$ type Dirichlet series defined as: \\

\begin{definition}
Let $\kappa>0,w \in \mathbb{C}, \alpha>0,\delta \geq 0,A \geq 0, B>0, M>0$ be some constants. A Dirichlet series $\mathcal{F}(s)$ defined as
$$ \mathcal{F}(s):=\sum_{n=1}^{\infty} f(n) n^{-s} $$
is said to be of type $\mathcal{P}(\kappa,w,\alpha,\delta,A,B,M)$ if the following conditions are satisfied:
\begin{enumerate}
\item for any $\epsilon >0$, we have
$$ \abs{f(n)} \ll_{\epsilon} n^{\epsilon} \qquad (n \geq 1);  $$

\item we have
$$ \sum_{n=1}^{\infty} \abs{f(n)}n^{-\sigma} \ll (\sigma-1)^{-\alpha} \qquad (\sigma>1); $$

\item the Dirichlet series
$$ \mathcal{G}(s;\kappa,w) := \mathcal{F}(s) \zeta(s)^{-\kappa} \zeta(2s)^w $$
is analytically continued to a holomorphic function in (some open set containing) $\Re(s) \geq \frac{1}{2}$ and, in this region $\mathcal{G}(s;\kappa,w)$ satisfies the bound 
$$ \abs{\mathcal{G}(s;\kappa,w)} \leq M \left( \, \abs{\tau}+1 \right)^{\max\{ \delta(1-\sigma),0 \}} \left( \log \left( \, \abs{\tau}+1 \right) \right)^A \qquad (s=\sigma+i\tau)$$
uniformly for $0<\kappa \leq B$ and \,$\abs{w} \leq B$. \\
\end{enumerate}
\end{definition}

\noindent
From Theorem II.5.1 of \cite{Gt}, the function
$$ Z(s;z):= \left\{ (s-1)\zeta(s) \right\}^z \qquad (z \in \mathbb{C}) $$
is holomorphic in the disc $\abs{s-1}<1$, and admits the Taylor series expansion 
$$ Z(s;z)=\sum_{j=0}^{\infty} \frac{\gamma_j(z)}{j!} (s-1)^j ,$$
where the $\gamma_j(z)'s$ are entire functions of $z$ and satisfy:
for all $B>0$ and $\epsilon>0$, the estimate
$$ \frac{\gamma_j(z)}{j!} \ll_{B,\epsilon} (1+\epsilon)^j \qquad (j \geq 0,\, \abs{z} \leq B). $$
Under our hypothesis, the function $\mathcal{G}(s;\kappa,w)\zeta(2s)^{-w}Z(s;\kappa)$ is holomorphic in the disc $\abs{s-1}<\frac{1}{2}$ and
$$ \abs{\mathcal{G}(s;\kappa,w)\zeta(2s)^{-w}Z(s;\kappa)} \ll_{A,B,\delta,\epsilon} M $$
for $\abs{s-1} \leq \frac{1}{2} -\epsilon$, $0<\kappa \leq B$ and $\abs{w} \leq B$. \\

\begin{theorem} 
Let $\kappa>0$, $w \in \mathbb{C}$, $\alpha > 0$, $\delta \geq 0$, $A \geq 0$, $B>0$, $M>0$ be some constants. Let $\eta_1 > 0$ be such that 
$$ \abs{\zeta(\sigma + it)} \ll \left( \, \abs{t}+2 \right)^{\eta_1(1-\sigma)}\log \left( \, \abs{t}+2 \right) \qquad \text{for} \;\; \frac{1}{2} \leq \sigma \leq 1+ \frac{1}{\log \left( \, \abs{t}+2 \right)} . $$ 
Suppose that
$$ \mathcal{F}(s):=\sum_{n=1}^{\infty} f(n)n^{-s}$$
is a Dirichlet series of type $\mathcal{P}(\kappa,w,\alpha,\delta,A,B,M)$. Then for any $\epsilon > 0$ and sufficiently large $x \geq x_0(\epsilon, \kappa, A)$, we have:
$$ \sum_{x<n\leq x+y} f(n) = y (\log x)^{\kappa-1} \left\{ \sum_{l=0}^N \frac{\lambda_l(\kappa,w)}{(\log x)^l} + O \left( R_N(x,y) \right)  \right\} $$ 
uniformly for 
$$ x \geq y \geq x^{\theta(\kappa,\delta)+\epsilon} \geq 2 , \; N \geq 0 , \; 0< \kappa \leq B , \; \abs{w} \leq B ,$$
where
$$ \lambda_l(\kappa,w) := \frac{g_l(\kappa,w)}{\Gamma(\kappa-l)} ,$$
$$ R_N(x,y) := \frac{y}{x} \sum_{l=1}^{N+1} \frac{l \abs{\lambda_{l-1}(\kappa,w)}}{(\log x)^l} + \frac{(a_1 N+1)^{N+1}}{x^{1/2}} + M \left\{ \left( \frac{a_1 N+1}{\log x} \right)^{N+1} + e^{-a_2 \frac{\log x}{\log \log x}} \right\} $$
for some constants $a_1,a_2 > 0$ and

\begin{equation*}
\theta(\kappa,\delta):=
\begin{cases}
\mbox{\large$\frac{5\delta+55\epsilon+7}{5\delta+5\epsilon+12}$} & \text{if} \; \; \kappa \leq \frac{12}{5\eta_1} ,\\ \\
\mbox{\large$\frac{\eta_1 \kappa+\delta-1+11\epsilon}{\eta_1 \kappa+\delta+\epsilon}$} & \text{if} \; \; \kappa > \frac{12}{5\eta_1} . 
\end{cases} 
\end{equation*} \\

\end{theorem} 

\begin{remark} \normalfont
This improves Theorem 1.1 of \cite{ZcJw} (See also \cite{Aas}). It is easy to check in either case (whether $\kappa \leq \frac{12}{5\eta_1}$ or $\kappa > \frac{12}{5\eta_1}$) that
$$ \theta(\kappa,\delta) < \frac{5\kappa+15\delta+21}{5\kappa+15\delta+36} $$
of \cite{ZcJw} for $\eta_1=\frac{1}{3}$.
Thus the above theorem is an improvement over the short interval length. The implied O--constant depends on various parameters like $A,B,\epsilon,\delta,\eta$ etc.  \\
\end{remark}

\begin{remark} \normalfont
$\eta_1=\frac{1}{3}$ follows from Hardy's estimate
$$ \abs{\zeta \left( \frac{1}{2}+it \right)} \ll \left( \, \abs{t}+2 \right)^{\frac{1}{6}} \log \left( \, \abs{t}+2 \right) .$$
In fact, one may even take the best--known value $\eta_1 < \frac{1}{3}$ from the work of Bourgain in \cite{Jb}, giving \\
$$ \abs{\zeta \left( \frac{1}{2}+it \right) } \ll \abs{t}^{\frac{13}{84}+\epsilon}. $$
\end{remark}

\begin{remark} \normalfont
If one assumes the zero density hypothesis for $\zeta(s)$, then we have
$$ N(\sigma,T) \ll T^{2(1-\sigma)}(\log T)^A .$$
Thus the above theorem holds with 

\begin{equation*}
\theta(\kappa,\delta):=
\begin{cases}
\mbox{\large$\frac{1+\delta+11\epsilon}{2+\delta+\epsilon}$} & \text{if} \; \; \kappa \leq \frac{2}{\eta_1}, \\ \\
\mbox{\large$\frac{\eta_1 \kappa+\delta-1+11\epsilon}{\eta_1 \kappa+\delta+\epsilon}$} & \text{if} \; \; \kappa > \frac{2}{\eta_1}.
\end{cases}
\end{equation*} \\

\end{remark}

\end{section}

\begin{section}{Construction of the Hooley--Huxley contour of integration} \hfill\\

\noindent
To construct the required Hooley--Huxley contour for our situation, we follow certain descriptions from H. Maier and A. Sankaranarayanan in \cite{HmAs}.\\
Let $C^*$ be a generic absolute constant in the following, which need not be the same at each occurrence.

\begin{definition}
\normalfont{
A zero $\rho = \beta +i \gamma $ $\left( \text{with} \ \beta \geq \frac{1}{2} \right)$ of $\zeta(s)$ is said to be \emph{good} if $ \beta < 1-\frac{C^*}{\log\log(\, \abs{t}+2)}$ and $\rho$ is said to be \emph{exceptional} otherwise. } \\
\end{definition}

\noindent
Let $T \geq T_0$ and $x \geq x_0$ ($T_0$ and $x_0$ are sufficiently large). Let $\mathcal{G}$ and $\mathcal{E}$ denote the set of all good and exceptional zeros of $\zeta(s)$ respectively with $\abs{\gamma} \leq T$. We denote by $|\mathcal{G}|$ and $|\mathcal{E}|$ to mean the cardinality of the sets $\mathcal{G}$ and $\mathcal{E}$ respectively. \\

\noindent
Let $\alpha$ be any fixed constant satisfying $\frac{1}{2}+\eta \leq \alpha \leq 1-\eta$ with $\eta$ being any arbitrarily small fixed positive constant. Since the contour will be symmetric with respect to the real axis, it suffices to describe it in the upper half--plane. We assume that $\abs{\mathcal{E}}=0$. Hence, $\zeta(s) \neq 0$ in the region $\left\{ \sigma > 1-\frac{C^{**}}{\log\log(U+12)} ,\; U \leq t \leq 2U \right\}$ where $C^{**}$ is a suitable absolute positive constant and we construct the contour accordingly. \\

\noindent
Let $T=2^{l_0}$. We choose $c$ with $\frac{1}{2} \leq c \leq 1$ such that $H_0=c\log\log T=2^L$ with a positive integer $L$. For $l \geq L$, write $U=U^{(l)}=2^l$. We define the contour for $U \leq t \leq 2U$. Let $H=H(U^{(l)})=c_l \log\log (U^{(l)})$ and choose $c_l$ satisfying $\frac{1}{2} \leq c_l \leq 1$ such that $\frac{U}{2H}$ is a positive integer. \\

\noindent
We split the interval $[U,\,2U]$ into $\frac{U}{2H}$ disjoint abutting small intervals $I_j = I_j^{(l)}$ of equal length $2H$ for $1 \leq j \leq \frac{U}{2H}$.
Let $I_j=[U_j-H,\,U_j+H]$ and let 
$$\beta_j = \sup \left\{ \beta \mid \rho = \beta +i \gamma ,\; \zeta(\rho)=0 ,\; \beta \geq \alpha ,\; \gamma \in [U_j-2H,\,U_j+2H] \, \right\}$$ 
and 
$$\beta_j^*=\beta_j+ \frac{C^*}{\log\log 2(U+12)}.$$

\noindent
We also define $\left( \text{with} \ H_0'=H_0+2(\log H_0)^2 \right)$
$$\beta_0 = \sup \left\{ \beta \mid \rho = \beta +i \gamma ,\; \zeta(\rho)=0 ,\; \beta \geq \alpha ,\; \gamma \in [0,\,2H_0'] \, \right\}$$ 
and
$$\beta_0^*=\beta_0+ \frac{C^*}{\log\log 2H_0}.$$

\noindent
If there is no zero of $\zeta(s)$ in the rectangle $\left\{ \sigma \geq \alpha,\; U_j-2H \leq t \leq U_j+2H \right\}$, then we define $\beta_j^*=\alpha$. A similar notion applies to $\beta_0^*$ too. \\

\noindent
Then the contour $\mathcal{C}$ consists of 
\begin{enumerate}
\item Vertical pieces ($V_j$):\\
\begin{equation*}
V_j=
\begin{cases}
[\beta_j^*+i(U_j-H+\epsilon),\; \beta_j^*+i(U_j+H-\epsilon)] & \text{if} \; \; \beta_j^* < \min{(\beta_{j-1}^*,\beta_{j+1}^*)} \\
[\beta_j^*+i(U_j-H-\epsilon),\; \beta_j^*+i(U_j+H+\epsilon)] & \text{if} \; \; \beta_j^* > \max{(\beta_{j-1}^*,\beta_{j+1}^*)} \\
[\beta_j^*+i(U_j-H-\epsilon),\; \beta_j^*+i(U_j+H-\epsilon)] & \text{if} \; \; \beta_{j-1}^* < \beta_j^* < \beta_{j+1}^* \\
[\beta_j^*+i(U_j-H+\epsilon),\; \beta_j^*+i(U_j+H+\epsilon)] & \text{if} \; \; \beta_{j+1}^* < \beta_j^* < \beta_{j-1}^* 
\end{cases}
\end{equation*}
and \\
\begin{equation*}
V_0=
\begin{cases}
[\beta_0^*,\; \beta_0^*+i(H_0-\epsilon)] & \text{if} \; \; \beta_1^* > \beta_0^* \\
[\beta_0^*,\; \beta_0^*+i(H_0+\epsilon)] & \text{if} \; \; \beta_1^* < \beta_0^* \\
\end{cases}
\end{equation*} \\

\item Horizontal pieces ($h_j$): \\
\begin{enumerate}
\item If $\beta_j^* < \min{(\beta_{j-1}^*,\beta_{j+1}^*)}$, then
\begin{align*}
h_j(\text{top}) &= [\beta_j^*+i(U_j+H-\epsilon),\; \beta_{j+1}^*+i(U_j+H-\epsilon)]  \quad \text{and} \\
h_j(\text{bottom}) &= [\beta_j^*+i(U_j-H+\epsilon),\; \beta_{j-1}^*+i(U_j-H+\epsilon)] \\
\end{align*}

\item If $\beta_j^* > \max{(\beta_{j-1}^*,\beta_{j+1}^*)}$, then
\begin{align*}
h_j(\text{top}) &= [\beta_{j+1}^*+i(U_j+H+\epsilon),\; \beta_j^*+i(U_j+H+\epsilon)]  \quad \text{and} \\
h_j(\text{bottom}) &= [\beta_{j-1}^*+i(U_j-H-\epsilon),\; \beta_j^*+i(U_j-H-\epsilon)] \\
\end{align*}

\item If $\beta_{j-1}^* < \beta_j^* < \beta_{j+1}^*$, then
\begin{align*}
h_j(\text{top}) &= [\beta_j^*+i(U_j+H-\epsilon),\; \beta_{j+1}^*+i(U_j+H-\epsilon)]  \quad \text{and} \\
h_j(\text{bottom}) &= [\beta_{j-1}^*+i(U_j-H-\epsilon),\; \beta_j^*+i(U_j-H-\epsilon)] \\
\end{align*}

\item If $\beta_{j+1}^* < \beta_j^* < \beta_{j-1}^*$, then
\begin{align*}
h_j(\text{top}) &= [\beta_{j+1}^*+i(U_j+H+\epsilon),\; \beta_j^*+i(U_j+H+\epsilon)]  \quad \text{and} \\
h_j(\text{bottom}) &= [\beta_j^*+i(U_j-H+\epsilon),\; \beta_{j-1}^*+i(U_j-H+\epsilon)] \\
\end{align*}

\end{enumerate}

\noindent
and similar horizontal pieces $h_{0,l}$ that link the top (respectively the bottom) vertical pieces of the ranges
$$ U^{(l-1)} \leq t \leq 2U^{(l-1)} \; \text{(respectively)} \; U^{(l)} \leq t \leq 2U^{(l)} .$$ \\
\end{enumerate}

\noindent
The vertical piece $V_0$ and the horizontal piece $h_0$ pertain to the interval $[2^3,H_0']$ where $H_0'=H_0+2(\log H_0)^2$. We also observe that the vertical piece $V^*$ for the interval $[0,2^3]$ can be taken to be $\alpha=\frac{1}{2}+\eta$ for any small positive constant $\eta$. Therefore, the contour $\mathcal{C}$ can be pictorially seen as shown below and
$$ \mathcal{C} = \Gamma \cup \Gamma_1 \cup \Gamma_2 \cup \Gamma_1^R \cup \Gamma_2^R .$$ 

\begin{figure}

\begin{tikzpicture}

\draw (2,-10) rectangle (13,10);
\draw[thick,->] (-1,0) -- (15,0) node[right]{$\sigma$};
\draw[thick,->] (0,-11) -- (0,11) node[above]{$\tau$};
\draw (0 pt,10 cm) node[anchor=east] {$T$};
\node at (0,10) {$\boldsymbol{\cdot}$};
\draw (0 pt,0 cm) node[anchor=north east] {$O$};
\node at (0,0) {$\boldsymbol{\cdot}$};
\draw (2 cm,0 pt) node[anchor=north east] {$\frac{1}{2}$};
\node at (2,0) {$\boldsymbol{\cdot}$};
\draw (13 cm,0 pt) node[anchor=north west] {$1+\frac{20}{\log x}$};
\node at (13,0) {$\boldsymbol{\cdot}$};
\draw (3 cm,3 pt) node[anchor=north] {$\frac{1}{2}+\eta$};
\node at (3,0) {$\boldsymbol{\cdot}$};
\draw (10 cm,0 pt) node[anchor=north] {$1$};
\node at (10,0) {$\boldsymbol{\cdot}$};
\draw (10.8 cm,0 pt) node[anchor=north west] {$1+\frac{1}{\log x}$};
\node at (10.85,0) {$\boldsymbol{\cdot}$};
\draw [blue] (3,0.5) -- node[above] {$\Gamma$} (9,0.5) node[pos=0.5]{\tikz \draw[- angle 90] (1 pt, 0) -- (-1 pt, 0);};
\draw [blue] (3,-0.5) -- (9,-0.5) node[pos=0.5]{\tikz \draw[- angle 90] (-1 pt, 0) -- (1 pt, 0);};
\draw [blue] (9,-0.5) arc (-150:150:1cm);
\draw [green] (3,0.5) -- node[right]{$\Gamma_1$} (3,1.5) node[pos=0.5]{\tikz \draw[- angle 90] (0,-1 pt) -- (0,1 pt);};
\draw (0 pt,1.5 cm) node[anchor=east] {$2^3$};
\node at (0,1.5) {$\boldsymbol{\cdot}$};
\draw [green] (3,1.5) -- (6,1.5) node[pos=0.5]{\tikz \draw[- angle 90] (-1 pt, 0) -- (1 pt, 0);};
\draw [red] (6,1.5) -- node[right] {$\Gamma_2$} (6,4) node[pos=0.5]{\tikz \draw[- angle 90] (0,-1 pt) -- (0,1 pt);};
\draw (0 pt,4 cm) node[anchor=east] {$H_0'$};
\node at (0,4) {$\boldsymbol{\cdot}$};
\draw [red] (6,4) -- (3,4);
\draw [red] (3,4) -- (3,4.5);
\draw [red] (3,4.5) -- (8,4.5);
\draw [red] (8,4.5) -- (8,5);
\draw [red] (8,5) -- (6,5);
\draw [red] (6,5) -- (6,5.5);
\draw [red] (6,5.5) -- (4,5.5);
\draw [red] (4,5.5) -- (4,6);
\draw [red] (4,6) -- (5,6);
\draw [red] (5,6) -- (5,6.5);
\draw [red] (5,6.5) -- (7,6.5);
\draw [red] (7,6.5) -- (7,7);
\draw [red] (7,7) -- (3.5,7);
\draw [red] (3.5,7) -- (3.5,7.5);
\draw [loosely dotted] (3.5,7.5) -- (3.5,10);

\draw [orange] (3,-0.5) -- node[right]{$\Gamma_1^R$} (3,-1.5) node[pos=0.5]{\tikz \draw[- angle 90] (0,-1 pt) -- (0,1 pt);};
\draw [orange] (3,-1.5) -- (6,-1.5) node[pos=0.5]{\tikz \draw[- angle 90] (1 pt, 0) -- (-1 pt, 0);};
\draw [purple] (6,-1.5) -- node[right] {$\Gamma_2^R$} (6,-4) node[pos=0.5]{\tikz \draw[- angle 90] (0,-1 pt) -- (0,1 pt);};

\draw [purple] (6,-4) -- (3,-4);
\draw [purple] (3,-4) -- (3,-4.5);
\draw [purple] (3,-4.5) -- (8,-4.5);
\draw [purple] (8,-4.5) -- (8,-5);
\draw [purple] (8,-5) -- (6,-5);
\draw [purple] (6,-5) -- (6,-5.5);
\draw [purple] (6,-5.5) -- (4,-5.5);
\draw [purple] (4,-5.5) -- (4,-6);
\draw [purple] (4,-6) -- (5,-6);
\draw [purple] (5,-6) -- (5,-6.5);
\draw [purple] (5,-6.5) -- (7,-6.5);
\draw [purple] (7,-6.5) -- (7,-7);
\draw [purple] (7,-7) -- (3.5,-7);
\draw [purple] (3.5,-7) -- (3.5,-7.5);
\draw [loosely dotted] (3.5,-7.5) -- (3.5,-10);

\end{tikzpicture}

\caption{Contour $\mathcal{C}$}
\end{figure}

\end{section}

\begin{section}{Proof of Theorem 1.1} \hfill\\

\begin{subsection}{Treatment of the sum \mbox{\boldmath$\sum \limits_{x<n \leq x+y} f(n)$}} \hfill\\

\noindent
Since $\mathcal{F}(s)$ is a Dirichlet series of the type $\mathcal{P}(\kappa,w, \alpha, \delta, A, B, M)$, we apply Corollary II.2.2.1 of \cite{Gt} with the choice of parameters $\sigma_a=1, B(n)=n^{\epsilon},\alpha=\alpha, \sigma=0$ to obtain

$$ \sum_{x<n \leq x+y} f(n) = \frac{1}{2\pi i} \int_{b-iT}^{b+iT} \mathcal{F}(s) \frac{(x+y)^s-x^s}{s} ds + O \left( \frac{x^{1+\epsilon}}{T} \right)$$
where $b=1+\frac{20}{\log x},100 \leq T \leq x$ such that $\zeta(\sigma+iT) \neq 0$ for $0<\sigma<1$.

\noindent
Now we replace the path of integration $[b-iT,b+iT]$ by the contour $\mathcal{C}$ described above. \\

\noindent
K. Ramachandra and A. Sankaranarayanan (see Theorems 1 and 2 of \cite{KrAs}) investigated certain upper bound estimates ``locally'' for the function \,$\abs{\log \mathcal{F}(s)}$ (where $\mathcal{F}(s)$ is any Dirichlet series satisfying certain general conditions) under the assumption that $\mathcal{F}(s) \neq 0$ in the rectangle $\{\sigma \geq \frac{1}{2}+\eta \; , \; T-H \leq t \leq T+H\}$ of $t$--width $2H$. Here the parameter $H$ can be chosen as small as $H=c \log \log \log T$. We record here a special case of the general theorem as: \\

\begin{lemma}
Let $\frac{1}{2} \leq \alpha^* \leq 1-\eta$, $H=a_3 \log \log T $ and suppose that $\zeta(s) \neq 0$ in $\{ \sigma > \alpha^* \; , \; \; T-H \leq t \leq T+H \}$. Then for $\alpha^* < \sigma \leq 1-\frac{a_4}{\log \log T} , \; \; T-\frac{H}{2} \leq t \leq T+\frac{H}{2} $, we have
$$ \abs{\log \zeta(\sigma+it)} \leq a_5 \log T (\log \log T)^{-1} $$ 
where $a_3$, $a_4$ and $a_5$ are certain positive constants.\\
\end{lemma}

\noindent
Therefore by this lemma, for $\abs{t}(\geq H_0)$, we have 
$$ \abs{\zeta(\sigma+it)} \ll U^{\epsilon} $$
for $\sigma+it \in V_j$ and $\sigma+it \in h_j$. The horizontal slab with $\abs{t} \in \left[ \frac{1}{\log x},\, H_0 \right]$ is treated as follows. We redefine 
 
\begin{equation*}
\beta_0^*=
\begin{cases}
\beta_{0,1}^* = \frac{1}{2}+\eta & \text{if} \; \; \abs{t} \in \left[ \frac{1}{\log x},\, 2^7 \right], \\
\beta_{0,2}^* = \beta_0 + \frac{C^*}{\log \log H_0'} & \text{if} \; \; \abs{t} \in [2^3,\, H_0'].
\end{cases}
\end{equation*} \\

\noindent
For the portion $\abs{t} \in [2^5,\, H_0]$, we first observe that the region 
$$ \left\{ 10 \geq \sigma \geq \beta_{0,2}^* \geq \beta_0+\frac{C^*}{\log \log H_0} \; , \; \; 2^3 \leq \abs{t} \leq H_0' \right\} $$
is free from zeros of $\zeta(s)$.
Therefore, applying the Borel--Carath\'eodory theorem, we get $\left( \text{for} \; 10 \geq \sigma \geq \beta_0+\frac{C^*}{\log \log H_0} \right),$
$$ \abs{\log \zeta(\sigma+it)} \ll (\log H_0)^{1-\epsilon} \ll (\log \log \log T)^{1-\epsilon} \ll H_0^{\epsilon}. $$
So this estimate holds when $\beta_{0,2}^*+it \in V_0$ with $t \in [2^5,H_0]$ and $T \geq T_0$ where $T_0$ is sufficiently large. \\

\noindent
The portion $\abs{t} \in [0,2^5]$ is dealt as follows. We observe that the region
$$ \left\{ 10 \geq \sigma \geq \beta_{0,1}^* \; , \; \; \abs{t} \leq 2^7 \right\} $$ 
is zero--free for $\zeta(s)$. This follows from the computational results. Thus,
$$ \abs{\zeta(s)} \ll \log \log (U+2) \ll U^{\epsilon} $$
$$ \text{for} \; \; \left\{ \sigma \geq \beta_{0,1}^* \; , \; \; \abs{s-1} \geq \frac{1}{10} \; , \; \; \abs{t} \leq 2^7 \right\} .$$
For $\abs{s-1} \leq \frac{1}{10}$,
$$ \abs{\zeta(s)} \ll \frac{1}{\abs{s-1}} \ll 1. $$ \\

\noindent
Thus, we need to estimate

\begin{align*}
\sum_{x<n \leq x+y} f(n) &= \frac{1}{2\pi i} \left[ \int_{\Gamma} + \int_{\Gamma_1} + \int_{\Gamma_2} + \int_{\Gamma_1^R} + \int_{\Gamma_2^R} \right] \mathcal{F}(s) \frac{(x+y)^s-x^s}{s} ds + O \left( \frac{x^{1+\epsilon}}{T} \right) \\
&= I_0+I_1+I_2+I_1^R+I_2^R + O \left( \frac{x^{1+\epsilon}}{T} \right) \; \text{(say)}. \\
\end{align*}

\end{subsection}

\begin{subsection}{Evaluation of \mbox{\boldmath$I_0$}} \hfill\\
\noindent
Let $0<a_6<\frac{1}{10}$ be any small constant. Since $\mathcal{G}(s;\kappa,w)\zeta(2s)^{-w}Z(s;\kappa)$ is holomorphic and $O(M)$ in the disc $\abs{s-1} \leq a_6$, the Cauchy's formula implies that
$$ g_l(\kappa,w) \ll M{a_6}^{-l} \qquad (l \geq 0, \; 0<\kappa \leq B, \; \abs{w}\leq B) $$
where $g_l(\kappa,w)$ is defined by
$$ \mathcal{G}(s;\kappa,w)\zeta(2s)^{-w}Z(s;\kappa) = \sum_{l=0}^{\infty} g_l(\kappa,w)(s-1)^l $$
with 
$$ Z(s;\kappa) := \left( (s-1)\zeta(s) \right)^{\kappa} ,$$
$$ g_l(\kappa,w) := \frac{1}{l!} \sum_{j=0}^l \binom lj \frac{\partial^{l-j} \left\{ \mathcal{G}(s;\kappa,w)\zeta(2s)^{-w} \right\} }{\partial s^{l-j}}\bigg|_{s=1} \gamma_j(\kappa).$$ \\

\noindent
Hence for any integer $N \geq 0$ and $\abs{s-1} \leq \frac{a_6}{2}$,
$$ \mathcal{G}(s;\kappa,w)\zeta(2s)^{-w}Z(s;\kappa) = \sum_{l=0}^N g_l(\kappa,w)(s-1)^l + O \left(M \left( \frac{\abs{s-1}}{a_6} \right)^{N+1} \right) .$$ \\
We have,
\begin{align*}
\mathcal{F}(s) &= \mathcal{G}(s;\kappa,w)\zeta(2s)^{-w} \zeta(s)^\kappa ,\\
&= \mathcal{G}(s;\kappa,w)\zeta(2s)^{-w}Z(s;\kappa)(s-1)^{-\kappa}.
\end{align*}

\noindent
Thus,
\begin{align*}
I_0 &:= \frac{1}{2 \pi i} \int_{\Gamma} \mathcal{F}(s) \frac{(x+y)^s-x^s}{s} ds \\
&= \sum_{l=0}^N g_l(\kappa,w) \frac{1}{2\pi i} \int_{\Gamma} (s-1)^{l-\kappa} \frac{(x+y)^s-x^s}{s} ds \\
&\quad + O \left(M{a_6}^{-N} \int_{\Gamma} (s-1)^{N+1-\kappa} \frac{(x+y)^s-x^s}{s} ds \right) \\
&= \sum_{l=0}^N g_l(\kappa,w) M_l(x,y) + O \left(M{a_6}^{-N} E_N(x,y) \right) \; \text{(say)}. \\
\end{align*}

\begin{subsubsection}{Evaluation of $M_l(x,y)$} \hfill\\
$$ M_l(x,y) := \frac{1}{2\pi i} \int_{\Gamma} (s-1)^{l-\kappa} \frac{(x+y)^s-x^s}{s} ds $$
Observe that
$$ \frac{(x+y)^s-x^s}{s} = \int_x^{x+y} u^{s-1} \ du .$$
Using Corollary II.5.2.1 of \cite{Gt}, we can write

\begin{align*}
M_l(x,y) &= \int_x^{x+y} \left( \frac{1}{2\pi i} \int_{\Gamma} (s-1)^{l-\kappa} u^{s-1} \ ds \right) du \\
&= \int_x^{x+y} (\log u)^{\kappa-1-l} \left\{ \frac{1}{\Gamma(\kappa-l)} + O \left( \frac{(a_7 l+1)^l}{u^{\frac{1}{2}}} \right) \right\} du
\end{align*}

\noindent
where we have used
$$ 47^{\, \abs{\kappa-l}} \Gamma \left( 1+\abs{\kappa-l} \right) \ll_B (a_7 l+1)^l \qquad (l \geq 0, \; 0<\kappa \leq B). $$
$a_7$ and the implied constant may depend at most on $B$.
Now for $0<\kappa \leq B$, $0<u<y \leq x$,

\begin{align*}
\log(x+u) &= \log x + \log \left( 1+\frac{u}{x} \right) \\
&= \log x+ O \left( \frac{u}{x} \right).
\end{align*}

\noindent
Therefore,
$$ \left( \log (x+u) \right)^{\kappa -1 -l} = (\log x)^{\kappa -1 -l} + O \left( \frac{(l+1)u(\log x)^{\kappa-2-l}}{x} \right) $$
and
\begin{align*}
\int_x^{x+y} (\log u)^{\kappa-1-l} du &= \int_0^y \left( \log (x+u) \right)^{\kappa-1-l} du \\
&= y(\log x)^{\kappa-1-l} +O \left( \frac{(l+1)(\log x)^{\kappa-2-l}}{x} \int_0^y u \ du \right) \\
&= y(\log x)^{\kappa-1-l} +O \left( \frac{(l+1)(\log x)^{\kappa-2-l}}{x} y^2 \right) \\
&= y(\log x)^{\kappa-1-l} \left\{1 + O_B \left( \frac{(l+1)y}{x \log x} \right) \right\} .
\end{align*}

\noindent
Also,
\begin{align*}
(a_7 l+1)^l \int_x^{x+y} \frac{(\log u)^{\kappa-1-l}}{u^{\frac{1}{2}}} du &\ll \frac{(a_7 l+1)^l}{x^{\frac{1}{2}}} \left( \log (2x) \right)^{\kappa-1-l}y \\
&\ll_B \frac{(a_7 l+1)^l (\log x)^{\kappa-1-l} y}{x^{\frac{1}{2}}}. \\
\end{align*}

\noindent
Thus, we get
$$ M_l(x,y) = y(\log x)^{\kappa-1-l} \left\{ \frac{1}{\Gamma(\kappa-l)} + O_B \left( \frac{(l+1)y}{\Gamma(\kappa-l)x \log x} \right) + O_B \left( \frac{(a_7 l+1)^l}{x^{\frac{1}{2}}} \right) \right\} $$
for $l \geq 0$, $0<\kappa \leq B$. \\

\end{subsubsection}

\begin{subsubsection}{Estimation of $E_N(x,y)$} \hfill\\

$$ E_N(x,y) := \int_{\Gamma} (s-1)^{N+1-\kappa} \frac{(x+y)^s-x^s}{s} ds $$

\noindent
We observe that 
\begin{align*}
\abs{\frac{(x+y)^s-x^s}{s}} &= \abs{\int_x^{x+y} u^{s-1} \ du} \leq \int_x^{x+y} u^{\sigma-1} \ du \\
&= \frac{u^{\sigma}}{\sigma}\bigg|_x^{x+y} = \frac{(x+y)^{\sigma}-x^{\sigma}}{\sigma} \\
&\ll \frac{x^{\sigma-1}y \sigma}{\sigma} \ll x^{\sigma-1}y. \\
\end{align*}

\noindent
Therefore, for $r=\frac{1}{\log x}$,
\begin{align*}
E_N(x,y) &\ll \int_{\frac{1}{2}+\eta}^{1-\frac{1}{\log x}} (1-\sigma)^{N+1-\kappa} x^{\sigma -1} y \ d \sigma \\
&\quad + \abs{ \int_{-\pi}^{\pi} (re^{i\theta})^{N+1-\kappa} \frac{(x+y)^{1+re^{i\theta}}-x^{1+re^{i\theta}}}{1+re^{i\theta}} re^{i\theta}i \ d\theta } \\
&\ll \frac{y}{(\log x)^{N+1-\kappa}} \int_{\frac{1}{2}}^{1-\frac{1}{\log x}} \left\{ (1-\sigma)\log x \right\}^{N+1-\kappa} e^{-(1-\sigma) \log x} \ d\sigma \\
&\quad + \int_{-\pi}^{\pi} \abs{r}^{N+1-\kappa} x^{r \cos \theta}y r \ d\theta  \\
&\ll \frac{y}{(\log x)^{N+1-\kappa}} \int_{1}^{\frac{\log x}{2}} u^{N+1-\kappa}e^{-u} \frac{du}{\log x} + yr^{N+2-\kappa} \\
&\ll \frac{y}{(\log x)^{N+2-\kappa}} \Gamma \left( 1+\abs{N-\kappa} \right) + \frac{y}{(\log x)^{N+2-\kappa}} \\
&\ll y(\log x)^{\kappa-1} \frac{(a_7 N+1)^{N+1}}{(\log x)^{N+1}}
\end{align*}

\noindent
uniformly for $x \geq y \geq 2 \; , \; N \geq 0 \; \text{and} \; 0<\kappa \leq B$ where $a_7>0$ and the implied constant depends only on $B$. 
Inserting all these estimates, we get

$$ I_0 = y(\log x)^{\kappa-1} \left\{ \sum_{l=0}^N \frac{\lambda_l(\kappa,w)}{(\log x)^l} + O_B \left( E_N^*(x,y) \right) \right\} $$
where 
$$ E_N^*(x,y) := \frac{y}{x} \sum_{l=1}^{N+1} \frac{l \abs{\lambda_{l-1}(\kappa,w)}}{(\log x)^l} + \frac{(a_7 N+1)^{N+1}}{x^{\frac{1}{2}}} + M\left( \frac{a_7 N+1}{\log x} \right)^{N+1}. $$ 
(This constant $a_7$ is denoted as $a_1$ in the statement of the theorem.) \\

\end{subsubsection}
\end{subsection}

\begin{subsection}{Treatment of \mbox{\boldmath$I_1$} and \mbox{\boldmath$I_1^R$}} \hfill\\
$$ I_1 := \frac{1}{2\pi i} \int_{\Gamma_1} \mathcal{F}(s) \frac{(x+y)^s-x^s}{s} ds $$
Note that
$$ \mathcal{F}(s) := \mathcal{G}(s;\kappa,w) \zeta(s)^\kappa \zeta(2s)^{-w} $$
so that
$$ \abs{\mathcal{G}(s;\kappa,w)} \leq M \left( \, \abs{\tau}+1 \right)^{\max \{ \delta(1-\sigma),0 \} } \left( \log \left( \, \abs{\tau}+1 \right) \right)^A \qquad (s=\sigma+i\tau),$$

\begin{align*}
\abs{\mathcal{G} \left( \frac{1}{2}+\eta+i\tau;\kappa,w \right) } &\leq M \left( \, \abs{\tau}+1 \right)^{\frac{\delta}{2}} \left( \log \left(\, \abs{\tau}+1 \right) \right)^A \\
&\leq M.2^{\delta2^3}(\log 2^4)^A \\
&\leq M.2^{8\delta}4^A \qquad \text{if} \; \abs{\tau} \leq 2^3. \\
\end{align*}

\noindent
In $\sigma>0,$ $\zeta(s)$ admits an analytic continuation as a single--valued function having its only singularity at $s=1$, which is a simple pole and one has the representation (in $\sigma>0$):

$$\zeta(s) = \frac{s}{s-1} - s \int_1^{\infty} \frac{(x)}{x^{s+1}} dx \ , \qquad  (x) \; \text{is the fractional part of} \; x $$ 
\begin{align*}
\abs{\zeta(\sigma+i\tau)} &\leq \frac{\sigma+\abs{\tau}}{\sqrt{(1-\sigma)^2+{\tau}^2}} + \left( \sigma+\abs{\tau} \right) \int_1^{\infty} \frac{dx}{x^{\sigma+1}} \\
&\leq \frac{\sigma+\abs{\tau}}{\sqrt{(1-\sigma)^2+{\tau}^2}} + \frac{\sigma+\abs{\tau}}{\sigma} \qquad \text{for} \; \sigma \geq \eta > 0. \\
\end{align*}

\noindent
Thus, 
$$ \abs{\zeta \left( \frac{1}{2}+\eta+i\tau \right)} \leq 2^6 $$ 
for $\abs{\tau} \leq 2^3$, $\eta>0$ be any small positive constant.
For $\kappa>0$,
$$ \abs{\zeta \left(\frac{1}{2}+\eta+i\tau \right)^{\kappa}} \leq 2^{6\kappa} .$$
For $w \in \mathbb{C}$,
$$\abs{\zeta(1+2\eta+2i\tau)^{-w}} \leq \abs{\zeta(1+2\eta+2i\tau)}^{a_8\, \abs{w}} \leq \left( \zeta(1+\eta) \right)^{a_8\, \abs{w}} $$
where $\abs{\tau} \leq 2^3$ and $a_8$ is an effective constant. \\

\noindent
Hence,
$$ \abs{\mathcal{F} \left( \frac{1}{2}+\eta+i\tau \right)} \leq M2^{8\delta}4^A 2^{6\kappa} \left( \zeta(1+\eta) \right)^{a_8\, \abs{w}} $$
and
\begin{align*}
\abs{I_1}+\abs{I_1^R} &\leq \frac{1}{2\pi} \abs{ \int_{-2^3}^{2^3} \mathcal{F} \left(\frac{1}{2}+\eta+i\tau \right) \frac{(x+y)^{\frac{1}{2}+\eta+i\tau}-x^{\frac{1}{2}+\eta+i\tau}}{\frac{1}{2}+\eta+i\tau} i \ d\tau} \\
&\leq \frac{1}{2\pi} M2^{8\delta}4^A2^{6\kappa} \left( \zeta(1+\eta) \right)^{a_8\, \abs{w}} \int_{-2^3}^{2^3} \frac{x^{\frac{1}{2}+\eta-1}y}{\abs{\frac{1}{2}+\eta+i\tau}} d\tau \\
&\leq M2^{8\delta}4^A 2^{6\kappa} \left( \zeta(1+\eta) \right)^{a_8\, \abs{w}} \frac{2^4}{\frac{1}{2}} \frac{y}{x^{\frac{1}{2}-\eta}} \\
&\ll_{A,B,\delta,\eta} M \frac{y}{x^{\frac{1}{2}-\eta}}
\end{align*}
uniformly for $0<\kappa \leq B$, $\abs{w}\leq B$. \\
\end{subsection}

\begin{subsection}{Estimation of the integral \mbox{\boldmath$I_2$}} \hfill\\
$$ I_2 := \frac{1}{2\pi i} \int_{\Gamma_2} \mathcal{F}(s) \frac{(x+y)^s-x^s}{s} ds $$

\noindent
Recall that
$$ \mathcal{F}(s) := \mathcal{G}(s;\kappa,w) \zeta(s)^{\kappa} \zeta(2s)^{-w} $$
in $\Gamma_2$, $\frac{1}{2}+\eta = \alpha \leq \sigma \leq 1-\frac{C^*}{\log \log U}$ for $t \in [U,\,2U]$ and $0<\kappa,\, \abs{w} \leq B$.
Also,
$$ \abs{\frac{(x+y)^s-x^s}{s}} \ll x^{\sigma-1}y, $$
$$ \abs{\mathcal{G}(s;\kappa,w)} \ll M \left( \, \abs{\tau}+1 \right)^{\max\{ \delta(1-\sigma),0 \}} \left( \log (\, \abs{\tau}+1) \right)^A \qquad \text{where} \; \;  s=\sigma+i\tau, $$
$$ \abs{\zeta(2s)^{-w}} \ll 1 .$$

\noindent
For $\frac{U}{2} \leq \abs{\tau} \leq 2U$, we have
$$ \abs{\zeta(\beta_j^*+i\tau)} \ll \log \log U $$ 
when $\beta_j^*$ is near to the left of the line $\sigma=1$ and
$$ \abs{\zeta(\beta_j^*+i\tau)} \ll e^{a_9 \frac{\log U}{\log \log U}} $$
when $\beta_j^*$ is away from the line $\sigma=1$ and closer to the line $\sigma=\frac{1}{2}$ from its right.

\noindent
The width of $V_j$ is $\ll H$, $\Re(V_j)=\beta_j^*$, $U_j-2H \leq \Im(V_j) \leq U_j+2H$, $H \ll a_3\log \log U$
and
$$ \mathcal{F}(s) \ll MU^{\delta(1-\beta_j^*)} \left( \log (U+1) \right)^A e^{a_9 \kappa \frac{\log U}{\log \log U}} .$$ \\

\noindent
The contribution of the vertical path $V_j$ to $I_2$ is
\begin{align*}
\abs{I_2(V_j)} &:= \abs{\frac{1}{2\pi i} \int_{V_j} \mathcal{F}(s) \frac{(x+y)^s-x^s}{s} ds} \\
&\ll MU^{\delta(1-\beta_j^*)} \left( \log (U+1) \right)^A e^{a_9 \kappa \frac{\log U}{\log \log U}} x^{\beta_j^*-1}y H \\
&\ll My \left( \frac{U^{\delta}}{x} \right)^{1-\beta_j^*} e^{a_9 \kappa \frac{\log U}{\log \log U}} \left( \log (U+1) \right)^A H \\
&\ll My \left( \frac{U^{\delta}}{x} \right)^{1-\beta_j^*} e^{a_9 \kappa \frac{\log U}{\log \log U}} (\log U)^{A+1} \\
&\ll My \ e^{a_{10} \frac{\log U}{\log \log U}} \left( \frac{U^{\delta}}{x} \right)^{1-\beta_j^*}.
\end{align*}

\noindent
Using $H_0 \ll U$, we get
\begin{align*}
\abs{I_2(V_0)} &:= \abs{\frac{1}{2\pi i} \int_{V_0} \mathcal{F}(s) \frac{(x+y)^s-x^s}{s} ds} \\
&\ll M {H_0}^{\delta(1-\beta_0^*)} \left( \log (H_0+1) \right)^A e^{a_9 \kappa \frac{\log U}{\log \log U}} x^{\beta_0^*-1}y H_0 \\
&\ll My \left( \frac{U^{\delta}}{x} \right)^{1-\beta_0^*} e^{a_9 \kappa \frac{\log U}{\log \log U}} \left( \log (U+1) \right)^A H_0 \\
&\ll My \left( \frac{U^{\delta}}{x} \right)^{1-\beta_0^*} e^{a_9 \kappa \frac{\log U}{\log \log U}} (\log U)^{A+1} \\
&\ll My \ e^{a_{10} \frac{\log U}{\log \log U}} \left( \frac{U^{\delta}}{x} \right)^{1-\beta_0^*}.
\end{align*}

\noindent
Thus, 
$$ \abs{I_2(V_j)} \ll_{A,B,\eta,\epsilon} My \ e^{a_{10} \frac{\log U}{\log \log U}} \left( \frac{U^{\delta}}{x} \right)^{1-\beta_j^*} \qquad \text{for} \; j=0,1,2,\dots, \frac{U}{2H}+1. $$\\

\noindent
Let $\beta^{**} = \max \{ \beta_j^*,\beta_{j+1}^* \}$. Then in $h_j$(top), $\alpha \leq \sigma \leq \beta^{**}$ and $\abs{\tau} = U_j+H - \epsilon \leq 2U$. 
The contribution of the horizontal path $h_j$ to $I_2$ is

\begin{align*}
\abs{I_2(h_j)} &:= \abs{\frac{1}{2\pi i} \int_{h_j} \mathcal{F}(s) \frac{(x+y)^s-x^s}{s} ds} \\
&\ll \int_{\alpha}^{\beta^{**}} MU^{\delta(1-\sigma)} \left( \log (U+1) \right)^A e^{a_9 \kappa \frac{\log U}{\log \log U}} x^{\sigma-1}y \ d\sigma \\
&\ll My \int_{\alpha}^{\beta^{**}} \left( \frac{U^{\delta}}{x} \right)^{1-\sigma} e^{a_9 \kappa \frac{\log U}{\log \log U}} \left( \log (U+1) \right)^A \ d\sigma \\
&\ll My \ e^{a_{10}\frac{\log U}{\log \log U}} \int_{\alpha}^{\beta^{**}} \left( \frac{U^{\delta}}{x} \right)^{1-\sigma} \ d\sigma .
\end{align*}

\noindent
Thus, 
$$ \abs{I_2(h_j)} \ll_{A,B,\eta,\epsilon} My \ e^{a_{10} \frac{\log U}{\log \log U}} \int_{\alpha}^{\beta^{**}} \left( \frac{U^{\delta}}{x} \right)^{1-\sigma} \ d\sigma \qquad \text{for} \; j=0,1,2,\dots, \frac{U}{2H}+1. $$

\noindent
Analogous estimate also applies for the horizontal pieces $h_{0,l}$. \\

\noindent
Let $\alpha$ be any fixed constant satisfying $\frac{1}{2}+\eta \leq \alpha \leq 1-\eta$ with $\eta$ being any arbitrarily small fixed positive constant. Assume that $|\mathcal{E}|=0$. We choose a partition of the interval $[\alpha,\,1]$ namely
$$\alpha=\alpha_0 < \alpha_1 < \dots < \alpha_{j-1} <\alpha_j=1 \qquad \text{with} \;  \alpha_j-\alpha_{j-1}<\epsilon.$$
The number of $j$--values for which $\beta_j^* \in [\alpha_{l-1},\,\alpha_l]$ is bounded by $N(\alpha_{l-1},2U)$.
Therefore, on the dyadic $t$--width $U \leq t \leq 2U$, the vertical bits and the horizontal bits contribute to the integral $I_2$, a quantity which is 

\begin{align*}
\abs{c(I_2)} &\leq \sum_j \left\{\, \abs{I_2(V_j)} + \abs{I_2(h_j)} \right\} \\
&\ll_{A,B,\eta,\epsilon} My \ e^{a_{10} \frac{\log U}{\log \log U}} \left[ \sideset{}{^*} \sum_j \left\{ \left( \frac{U^{\delta}}{x} \right)^{1-\beta_j^*} + \int_{\alpha}^{\beta^{**}} \left( \frac{U^{\delta}}{x} \right)^{1-\sigma} \ d\sigma \right\} \right] \\
&\ll_{A,B,\eta,\epsilon} My \ e^{a_{10} \frac{\log U}{\log \log U}} \int_{\alpha}^{1-\sigma_0} \left( \frac{U^{\delta}}{x} \right)^{1-\sigma} \ dN^*(\sigma,2U)
\end{align*}
where $N^*(\sigma,U) := \sideset{}{^*} \sum \limits_{\sigma \leq \beta_j^*, \atop \abs{\gamma_j} \leq U} 1$. 

\noindent
Thus,
\begin{align*} 
\abs{c(I_2)} &\ll_{A,B,\eta,\epsilon} My \ e^{a_{10} \frac{\log U}{\log \log U}} \left\{ \left( \frac{U^{\delta}}{x} \right)^{1-\sigma} N^*(\sigma,2U) \bigg|_{\alpha}^{1-\sigma_0} \right. \\
&\quad+ \left. \int_{\alpha}^{1-\sigma_0} \left( \frac{U^{\delta}}{x} \right)^{1-\sigma} \abs{\log \left( \frac{U^{\delta}}{x} \right) } N^*(\sigma,2U) \ d\sigma \right\} \\
&\ll_{A,B,\eta,\epsilon,\delta} My \ e^{a_{10} \frac{\log U}{\log \log U}} \left\{ \left( \frac{U^{\delta}}{x} \right)^{1-\alpha} N(\alpha,2U) \right. \\ 
&\quad+ \left. \log T \int_{\alpha}^{1-\sigma_0} \left( \frac{U^{\delta}}{x} \right)^{1-\sigma} N(\sigma,2U) \ d\sigma \right\} .\\
\end{align*}

\noindent
Here, $\sigma_0=\sigma_0(T) := \frac{C^*}{\log \log T}$.
According to our assumption $|\mathcal{E}|=0$, i.e., for $\sigma \geq 1-\sigma_0$, $\zeta(s) \neq 0$.
Therefore, $N^*(\sigma,2U) \leq N(\sigma,2U)$.
From \cite{Mn}, it is known that
$$ N(\sigma,T) \ll T^{\frac{12}{5}(1-\sigma)} (\log T)^{44} $$
for $\frac{1}{2} \leq \sigma \leq 1$ and $T \geq 2$.
Hence,
\begin{align*}
c(I_2) &\ll My \ e^{a_{10} \frac{\log U}{\log \log U}} \left\{ \left( \frac{U^{\delta}}{x} \right)^{1-\alpha} U^{\frac{12}{5}(1-\alpha)} (\log T)^{44} \right. \\
&\quad + \left. \log T \int_{\alpha}^{1-\sigma_0} \left( \frac{U^{\delta}}{x} \right)^{1-\sigma} U^{\frac{12}{5}(1-\sigma)} (\log T)^{44} \ d\sigma \right\} \\
&\ll My \ e^{a_{10} \frac{\log U}{\log \log U}} \left\{ \left( \frac{U^{\delta}}{x} \right)^{1-\alpha} U^{\frac{12}{5}(1-\alpha)} (\log T)^{44} \right. \\
&\quad+ \left. (\log T)^{46} \int_{\alpha}^{1-\sigma_0} \left( \frac{U^{\delta+\frac{12}{5}}}{x} \right)^{1-\sigma} \ d\sigma \right\} . \\
\end{align*}

\noindent
Note that as a function of $\sigma$, $\left( \frac{U^{\delta+\frac{12}{5}}}{x} \right)^{1-\sigma}$ is monotonic in $[\alpha,\,1-\frac{\sigma_0}{2}]$ and hence it attains maximum at its extremities. Therefore,

\begin{align*}
c(I_2) &\ll My \ e^{a_{10} \frac{\log U}{\log \log U}} \left\{ \left( \frac{U^{\delta+\frac{12}{5}}}{x} \right)^{1-\alpha} (\log T)^{46} + \left( \frac{U^{\delta+\frac{12}{5}}}{x} \right)^{\sigma_0} (\log T)^{46} \right\} \\
&\ll My \ e^{a_{10} \frac{\log U}{\log \log U}} (\log T)^{46} \left\{ \left( \frac{U^{\delta+\frac{12}{5}}}{x} \right)^{1-\alpha} + \left( \frac{U^{\delta+\frac{12}{5}}}{x} \right)^{\sigma_0} \right\} 
\end{align*}
and
$$ I_2 \ll \sum_{U=2^l, \atop \frac{\log T}{\log 2} \geq l \geq L} c(I_2) \ll My \ e^{a_{11} \frac{\log T}{\log \log T}} \left\{ \left( \frac{T^{\delta+\epsilon+\frac{12}{5}}}{x} \right)^{1-\alpha} + \left( \frac{T^{\delta+\epsilon+\frac{12}{5}}}{x} \right)^{\sigma_0} \right\} $$
since $T \ll U \ll T$. For the sake of convenience, we have multiplied the first term in the curly bracket by $T^{\epsilon (1-\alpha)}$ and the second term by $T^{\epsilon \sigma_0}$.

\noindent
A similar estimate holds for $I_2^R$, of course on the assumption that $|\mathcal{E}|=0$. \\

\end{subsection}

\begin{subsection}{Case \mbox{\boldmath$|\mathcal{E}| \geq 1$}} \hfill\\
\noindent
One observes that the number of exceptional zeros is
\begin{align*}
|\mathcal{E}| &:= N(1-\sigma_0,T) \\
&\ll T^{\frac{12}{5}\sigma_0} (\log T)^{44} \\
&\ll e^{\frac{C^* \log T}{\log \log T}} \\
&\ll T^{\epsilon}
\end{align*}
for $T \geq T_0$ ($T_0$ sufficiently large). \\

\noindent
Recall that
$$ \mathcal{F}(s) := \mathcal{G}(s;\kappa,w) \zeta(s)^{\kappa} \zeta(2s)^{-w} , $$
$$ \frac{1}{2}+\eta = \alpha \leq \sigma \leq 1-\frac{C^*}{\log \log U} \leq \beta_{j,e}^* ,$$
$$ \mathcal{G}(s;\kappa,w) \ll M \left( \, \abs{\tau}+1 \right)^{\max\{ \delta(1-\sigma),0 \}} \left( \log \left(\, \abs{\tau}+1\right) \right)^A \qquad \text{where} \; s=\sigma+i\tau \; \text{and} $$
$$ \abs{\zeta(s)}^\kappa \ll \left( \, \abs{\tau}+1 \right)^{\eta_1 \kappa(1-\sigma)} \log \left( \, \abs{\tau}+1 \right) \qquad \text{with} \; \eta_1 < \frac{1}{3} .$$ \\

\noindent
The contribution of the vertical path $V_{j,e}$ pertaining to an exceptional zero $\beta_{j,e}$ is
\begin{align*}
\abs{I_2(V_{j,e})} &\ll MU^{\delta(1-\beta_{j,e}^*)} \left( \log(U+1) \right)^A x^{\beta_{j,e}^*-1}y H U^{\eta_1 \kappa(1-\beta_{j,e}^*)} \log(U+1)\\
&\ll M U^{(\eta_1 \kappa +\delta)(1-\beta_{j,e}^*)} \left( \log(U+1) \right)^{A+1} y H x^{\beta_{j,e}^*-1} \\
&\ll M U^{(\eta_1 \kappa +\delta)(1-\beta_{j,e}^*)} \left( \log(U+1) \right)^{A+2} y x^{\beta_{j,e}^*-1} .
\end{align*}

\noindent
Similarly,
$$ \abs{I_2(V_{0,e})} \ll M H_0^{(\eta_1 \kappa +\delta)(1-\beta_{0,e}^*)} (\log T)^{A+2} x^{\beta_{0,e}^*-1} y. $$

\noindent
Therefore,
$$ \abs{I_2(V_{j,e})} \ll M U^{(\eta_1 \kappa +\delta)(1-\beta_{j,e}^*)} (\log T)^{A+2} y x^{\beta_{j,e}^*-1} $$
for $j=0,1,2,\dots,\frac{U}{2H}+1$. \\

\noindent
Similarly, the horizontal path $h_{j,e}$ contributes to $I_2$,
\begin{align*}
\abs{I_2(h_{j,e})} &\ll My \int_{\alpha}^{\beta_{j,e}^*} U^{(\eta_1 \kappa+\delta)(1-\sigma)} (\log T)^{A+2} x^{\sigma-1} \ d\sigma \\
&\ll My \int_{\alpha}^{\beta_{j,e}^*} U^{(\eta_1 \kappa+\delta)(1-\sigma)} (\log T)^{A+2} x^{\sigma-1} \ d\sigma 
\end{align*}
with $\alpha \leq \sigma \leq \beta_{j,e}^*$ and $1- \frac{C^*}{\log \log U} \leq \beta_{j,e}^*$. \\

\noindent
Thus in the case of exceptional set $\mathcal{E}$ being non--empty, we obtain as before
\begin{align*}
c_e(I_2) &\ll My(\log T)^{A+2} \left\{ \left( \frac{U^{\eta_1 \kappa +\delta}}{x} \right)^{1-\alpha} N_e^*(\alpha,2U) \right. \\
&\quad + \left. (\log T) \int_{\alpha}^1 \left( \frac{U^{\eta_1 \kappa +\delta}}{x} \right)^{1-\sigma} N_e^*(\sigma,2U) \ d\sigma \right\} \\
&\ll My(\log T)^{A+3} \left\{ \left( \frac{U^{\eta_1 \kappa+\delta+\epsilon}}{x} \right)^{1-\alpha} + \left( \frac{U^{\eta_1 \kappa+\delta+\epsilon}}{x} \right) +1 \right\} 
\end{align*}
where $N_e^*(\sigma,U) := {\sum \limits_{\sigma \leq \beta_{j,e}^*, \atop \abs{\gamma_{j,e}} \leq U} 1}$. Thus, $N_e^*(\sigma,U) \leq |\mathcal{E}| \ll U^{\epsilon}$. \\
Therefore, 
$$ I_{2,e} \ll \sum_{U=2^l,\atop \frac{\log T}{\log 2} \geq l \geq L} c_e(I_2) \ll My (\log T)^{A+4} \left\{ \left( \frac{T^{\eta_1 \kappa+\delta+\epsilon}}{x} \right)^{1-\alpha} + \left( \frac{T^{\eta_1 \kappa+\delta+\epsilon}}{x} \right) \right\} .$$

\noindent
A similar estimate also holds for $I_{2,e}^R$. \\ \\
We observe that
$$ N^*(\sigma,2U) = N_{\mathcal{G}}^*(\sigma,2U)+N_{\mathcal{E}}^*(\sigma,2U). $$
Thus we get in any case whether $|\mathcal{E}|=0$ or $|\mathcal{E}| \geq 1$,
\begin{align*}
I_2 &\ll My \ e^{a_{11} \frac{\log T}{\log \log T}} \left\{ \left( \frac{T^{\delta+\epsilon+\frac{12}{5}}}{x} \right)^{1-\alpha} + \left( \frac{T^{\delta+\epsilon+\frac{12}{5}}}{x} \right)^{\sigma_0} \right\} \\
&\quad+ My(\log T)^{A+4} \left\{ \left( \frac{T^{\eta_1 \kappa+\delta+\epsilon}}{x} \right)^{1-\alpha} + \left( \frac{T^{\eta_1 \kappa+\delta+\epsilon}}{x} \right) \right\} . \\
\end{align*}

\noindent
A similar estimate holds for $I_2^R$. Thus, 
\begin{align*}
I_1+I_1^R+I_2+I_2^R &\ll_{A,B,\delta,\eta} M\frac{y}{x^{\frac{1}{2}-\eta}} + My \ e^{a_{11} \frac{\log T}{\log \log T}} \left\{ \left( \frac{T^{\delta+\epsilon+\frac{12}{5}}}{x} \right)^{1-\alpha} \right. \\
&\quad + \left. \left( \frac{T^{\eta_1 \kappa+\delta+\epsilon}}{x} \right)^{1-\alpha} +  \left( \frac{T^{\delta+\epsilon+\frac{12}{5}}}{x} \right)^{\sigma_0}  + \left( \frac{T^{\eta_1 \kappa+\delta+\epsilon}}{x} \right) \right\} . \\
\end{align*} 
\end{subsection}

\begin{subsection}{Case 1: If \mbox{\boldmath $\kappa \leq \frac{12}{5\eta_1}$}}
\begin{align*}
I_1+I_1^R+I_2+I_2^R &\ll M\frac{y}{x^{\frac{1}{2}-\eta}} + My \ e^{a_{11} \frac{\log T}{\log \log T}} \left\{ \left( \frac{T^{\delta+\epsilon+\frac{12}{5}}}{x} \right)^{1-\alpha} \right. \\
&\quad + \left. \left( \frac{T^{\delta+\epsilon+\frac{12}{5}}}{x} \right)^{\sigma_0}  + \left( \frac{T^{\delta+\epsilon+\frac{12}{5}}}{x} \right) \right\}  \\
\end{align*}

\noindent
We choose $T$ such that $T^{\delta+\epsilon+\frac{12}{5}} \sim x^{1-10\epsilon}$ so that
\begin{align*}
I_1+I_1^R+I_2+I_2^R &\ll M\frac{y}{x^{\frac{1}{2}-\eta}} + My \ e^{a_{12} \frac{\log x}{\log \log x}} \left\{ x^{-10\epsilon(1-\alpha)}+x^{-10\epsilon \sigma_0}+x^{-10\epsilon} \right\} \\
&\ll My \ e^{-a_{13}(\alpha, \sigma_0, \epsilon) \frac{\log x}{\log \log x}}
\end{align*}
for sufficiently large $x \geq x_0(\epsilon, \kappa, A)$. 

\noindent
From the error term in the Perron's formula,
\begin{align*}
\frac{x^{\mbox{\normalfont$1+\epsilon$}}}{T} &\ll x^{\mbox{\normalfont $1+\epsilon - \frac{1-10\epsilon}{\delta+\epsilon+\frac{12}{5}}$}} \ll x^{\mbox{\normalfont $\theta +\epsilon$}} \\
x^{\mbox{\normalfont $\theta$}} &\gg x^{\mbox{\normalfont$ 1- \frac{1-10\epsilon}{\delta+\epsilon+\frac{12}{5}}$}} \\
\theta &\geq 1-\frac{5-50\epsilon}{5\delta+5\epsilon+12} \\
\theta &\geq \frac{5\delta+55\epsilon+7}{5\delta+5\epsilon+12} . \\
\end{align*}
\end{subsection}

\begin{subsection}{Case 2: If \mbox{\boldmath $\kappa > \frac{12}{5\eta_1}$}}
\begin{align*}
I_1+I_1^R+I_2+I_2^R &\ll M\frac{y}{x^{\frac{1}{2}-\eta}} + My \ e^{a_{11} \frac{\log T}{\log \log T}} \left\{ \left( \frac{T^{\eta_1 \kappa+\delta+\epsilon}}{x} \right)^{1-\alpha} \right. \\
&\quad + \left. \left( \frac{T^{\eta_1 \kappa+\delta+\epsilon}}{x} \right)^{\sigma_0}  + \left( \frac{T^{\eta_1 \kappa+\delta+\epsilon}}{x} \right) \right\}  \\
\end{align*}

\noindent
We choose $T$ such that $T^{\eta_1 \kappa+\delta+\epsilon} \sim x^{1-10\epsilon}$ so that
\begin{align*}
I_1+I_1^R+I_2+I_2^R &\ll M\frac{y}{x^{\frac{1}{2}-\eta}} + My \ e^{a_{14} \frac{\log x}{\log \log x}} \left\{ x^{-5\epsilon}+x^{-10\epsilon \sigma_0}+x^{-10\epsilon} \right\} \\
&\ll My \ e^{-a_{15}(\epsilon, \sigma_0) \frac{\log x}{\log \log x}} 
\end{align*}
for sufficiently large $x \geq x_0(\epsilon, \kappa, A)$. \\

\noindent
From the error term in the Perron's formula,
\begin{align*}
\frac{x^{\mbox{\normalfont$1+\epsilon$}}}{T} &\ll x^{\mbox{\normalfont $1+\epsilon - \frac{1-10\epsilon}{\eta_1 \kappa+\delta+\epsilon}$}} \ll x^{\mbox{\normalfont $\theta +\epsilon$}} \\
x^{\mbox{\normalfont $\theta$}} &\gg x^{\mbox{\normalfont $1-\frac{1-10\epsilon}{\eta_1 \kappa+\delta+\epsilon}$}} \\
\theta &\geq \frac{\eta_1 \kappa+\delta-1+11\epsilon}{\eta_1 \kappa+\delta+\epsilon} . 
\end{align*}

\end{subsection}

\end{section}

\noindent
This completes the proof of Theorem 1.1. \\

\begin{section}{Consequences of Hal\'asz--Tur\'an theorem} \hfill\\

\begin{theorem}\cite{GhPt}
Assume the Lindel\"{o}f hypothesis for $\zeta(s)$ in the form 
$$ \abs{\zeta \left( \frac{1}{2}+it \right)} \leq t^{\eta_2^2} \qquad \text{for} \; \;  t>t_0 $$
for all sufficiently small positive numbers $\eta_2$. Then the inequality
$$ N \left( \frac{3}{4}+2\eta_2^{\frac{1}{2}},T \right) < T^{3\eta_2} $$
holds for $T>t_0$.
\end{theorem}

\noindent
For a more general theorem along the same flavour see Theorem 1.3 and Theorem 1.4 of \cite{As}.\\

\noindent
Therefore, if one assumes Lindel\"{o}f hypothesis in the form given in the above theorem, then by taking $\alpha=\frac{3}{4} + 2\eta_2^{\frac{1}{2}}$ in the earlier arguments, we find that 
\begin{align*}
&I_1+I_1^R+I_2+I_2^R \ll_{A,B,\delta,\eta} M\frac{y}{x^{1-\alpha}} \\
&\quad +My \ (\log T)^{A+10} \left\{ \left( \frac{T^{2\eta_2^2 \kappa+\delta+\frac{3\eta_2}{(1-\beta^*)}}}{x} \right) ^{1-\alpha}+ \left( \frac{T^{2\eta_2^2 \kappa+\delta+\frac{3\eta_2}{(1-\beta^*)}}}{x} \right) \right\} . \\
\end{align*}

\noindent
Since for $\sigma \geq \frac{3}{4}+2\eta_2^{\frac{1}{2}}$, under Lindel\"{o}f hypothesis,
\begin{align*}
N(\sigma,T) &\leq N \left(\frac{3}{4}+2\eta_2^{\frac{1}{2}},T \right) \\
&< T^{3\eta_2} \\
&= T^{\frac{3\eta_2(1-\sigma)}{(1-\sigma)}} \\
&< T^{\frac{3\eta_2(1-\sigma)}{(1-\beta^*)}} .
\end{align*}
Here $\beta^*$ is defined to be slightly to the right of $\max \limits_j \beta_j^*$ and this is possible because of the standard zero-free region.\\

\noindent
From convexity principle for $\frac{1}{2} \leq \sigma \leq 1$, we find that
\begin{align*}
\zeta(\sigma+it) &\ll t^{\eta_2^2} x^{\frac{1}{2}-\sigma}+(\log t)x^{1-\sigma} \\
&\ll t^{2\eta_2^2(1-\sigma)}(\log t)
\end{align*}
by choosing $x=t^{2\eta_2^2}$. \\

\noindent
Now we choose $\eta_2=\epsilon(1-\beta^*)$ so that $\eta_2^2<\epsilon^2<\epsilon$ and

\begin{align*}
&I_1+I_1^R+I_2+I_2^R \ll M\frac{y}{x^{\frac{1}{4}-2\eta_2^{\frac{1}{2}}}} \\
&\quad +My \ (\log T)^{A+10} \left\{ \left( \frac{T^{2 \kappa \epsilon+\delta+3\epsilon}}{x} \right) ^{1-\alpha}+ \left( \frac{T^{2 \kappa \epsilon+\delta+3\epsilon}}{x} \right) \right\}.
\end{align*}

\noindent
By choosing $T^{2 \kappa \epsilon+\delta+3\epsilon} \sim x^{1-10\epsilon}$, we observe that

$$ I_1+I_1^R+I_2+I_2^R \ll M\frac{y}{x^{\frac{1}{4}-2\eta_2^{\frac{1}{2}}}}+My \ (\log x)^{A+10} \left\{ x^{-10\epsilon(1-\alpha)} + x^{-10\epsilon} \right\} .$$ \\

\noindent
Hence under the assumption of Lindel\"{o}f hypothesis, the above theorem holds with

$$ \theta(\kappa,\delta) := \frac{\delta-1+2\kappa \epsilon+13\epsilon}{\delta +2\kappa \epsilon +3\epsilon} .$$

\noindent
One observes that $\frac{\delta-1}{\delta} < \frac{1}{2}$ when $1<\delta<2$ and $\frac{\delta-1}{\delta}<1$ for any positive $\delta$. One needs to assume that $\delta>1$ so that the numerator is positive. \\ \\

\noindent
Applying the above contour with $\alpha=\frac{3}{4}+2\eta_2^{\frac{1}{2}}$ (assuming Lindel\"{o}f hypothesis in the stated form), we obtain

\begin{align*}
I_2 &\ll \sum_{\frac{\log T}{\log 2} \geq l \geq L} c(I_2) \\
& \ll My \ (\log T)^{A+3} \left\{ \left( \frac{T^{2\eta_2^2 \kappa+\delta}}{x} \right)^{1-\alpha} N^{**}(\alpha,T) \right. \\
&\quad + \left. (\log T) \int_{\alpha}^1 \left( \frac{T^{2\eta_2^2 \kappa+\delta}}{x} \right)^{1-\sigma} N^{**}(\sigma,T) \ d\sigma\right\} \\
&\ll My \ (\log T)^{A+10} T^{3\eta_2} \left\{ \left( \frac{T^{2\eta_2^2 \kappa+\delta}}{x} \right)^{1-\alpha} + \left( \frac{T^{2\eta_2^2 \kappa+\delta}}{x} \right) \right\} \\
&\ll My \ T^{3\eta_2+\epsilon} \left\{ \left( \frac{T^{2\eta_2^2 \kappa+\delta}}{x} \right)^{1-\alpha} + \left( \frac{T^{2\eta_2^2 \kappa+\delta}}{x} \right) \right\} . \\
\end{align*}
Here $N^{**}(\sigma,T)$ has its relevant meaning with the current context of $\alpha$.

\noindent
Choose $\eta_2=\epsilon$ and $T$ such that $T^{2\eta_2^2 \kappa+\delta} = x^{1-20\epsilon}$. Then,
\begin{align*}
I_2 &\ll My \ T^{4\epsilon} \{ x^{-20\epsilon(1-\alpha)}+x^{-20\epsilon} \} \\
&\ll My \ x^{4\epsilon} \{ x^{-5\epsilon+40\epsilon^{\frac{3}{2}}}+x^{-20\epsilon} \} . \\
\end{align*}

\noindent
Hence, we can also take
\begin{align*}
\theta(\kappa,\delta) &:=1-\frac{1-20\epsilon}{2\eta_2^2 \kappa+\delta} \\
&= \frac{2\eta_2^2 \kappa+\delta -1 +20\epsilon}{2\eta_2^2 \kappa+\delta} \\
&= \frac{2\epsilon^2 \kappa+\delta -1 +20\epsilon}{2\epsilon^2 \kappa+\delta}. \\
\end{align*}

\noindent
Again of course, one needs to assume that $\delta >1$.
We observe that the earlier unconditional estimate for $\theta(\kappa,\delta) \left(\text{relevant when} \ \kappa < \frac{12}{5\eta_1} \right)$ is $\frac{5\delta+7}{5\delta+12}+\epsilon_1$ which may be compared with the Lindel\"{o}f hypothesis conditional estimate $\frac{\delta-1}{\delta}+\epsilon_2$. Clearly,
$$ \frac{5\delta+7}{5\delta+12} > \frac{\delta-1}{\delta} $$
for any $\delta>0$. (However, our relevance here is $\delta>1$.) \\

\end{section}

\subsection*{Acknowledgements}
\noindent
The first author is grateful to UGC for its supporting NET fellowship with UGC Ref. No.: 1004/(CSIR--UGC NET Dec. 2017). The authors are thankful to the anonymous referee for some fruitful comments and for pointing out some discrepancies in the earlier versions. \\

\end{document}